\documentclass[11pt]{article}
\usepackage{a4wide}
\usepackage{amsmath}

\usepackage{amsfonts}
\usepackage{amssymb}
\usepackage{euscript}
\newenvironment{proof}{\noindent {\it Proof.~~}\ }{\  \rule{1mm}{2mm}\medskip}

\newtheorem{theorem}{Theorem}
\newtheorem{lemma}[theorem]{Lemma}
\newtheorem{corollary}[theorem]{Corollary}
\newtheorem{proposition}[theorem]{Proposition}

\newtheorem{theirtheorem}{Theorem}

\newtheorem{theirlemma}[theirtheorem]{Lemma}

\newcommand{\subgp}[1]{\langle{#1}\rangle}

\begin{document}
\title{
  Hyper-atoms applied to  the critical pair Theory}

\author{ Yahya O. Hamidoune\thanks{UPMC Univ Paris 06,
 E. Combinatoire, Case 189, 4 Place Jussieu,
75005 Paris, France,     {\tt hamidoune@math.jussieu.fr} }
}
\maketitle

\begin{abstract}

We introduce the notion of a hyper-atom and prove a basic property of this object.
This new method allows to improve several results in  the classical critical pair theory including its cornerstone:
the  Kemperman Structure Theorem.
\end{abstract}

\section{Introduction}

Let $G$ be an abelian group and
let  $ A$ and $B$ be subsets of $ G $. The subgroup generated by  $A$ will
be denoted by $\subgp{A}$. The {\em sumset} $A+B$ is defined as
$$A+B=\{x+y \ : \ x\in A\  \mbox{and}\ y\in
  B\}.$$
Let $H$ be a subgroup of $G.$
We shall say that $H$ is a {\em proper} subgroup if $H\neq G$.
We shall denote by  $\phi _H$  the canonical morphism from $G$ onto $G/H$.
We shall say that $A$ is $H$-periodic if
$A+H=A.$
The {\em period} of $A$ is
$G_A=\{x\in G :A+x=A\}$. A set having a non-zero period is said to be  {\em periodic}.
A non-periodic set is said to be  {\em aperiodic}.
A basic tool in Additive Number Theory is the following
generalization of the Cauchy-Davenport Theorem
due to Kneser:

\begin{theirtheorem}[Kneser \cite{tv}]\label{kneser}
Let  $A, B\subset G$ be finite
subsets of  an abelian group. If $A+B$ is aperiodic,  then $|A+B|\ge |A|+|B|-1$.

\end{theirtheorem}

The description of the subsets $A$ and $B$ with $|A+B|= |A|+|B|-1$, obtained by Kemperman in \cite{kempacta},  is a deep result in the classical critical pair theory. Another step in this direction is proposed by Grynkiewicz in \cite{davkem}. These two results are proved within about 80 pages.
One of our  aims in the present work is to present a methodology leading to generalizations, new  results and relatively short proofs.
The present   work is essentially self-contained. We assume only Kneser's Theorem, Theorem \ref{inter2frag} and Theorem \ref{2atomejc}. The last two results  are proved in around  2 pages  in \cite{hiso2007}.

Let $S$ be a proper subset of an abelian group $G$ with $0\in S$ and put ${\cal S}_k=\{X : |X|\ge k \ \mbox{and} \  |X+S|\leq |G|-k\}.$
We shall write  $$\kappa _k=
\min \{|X+S|-|X|: X\in {\cal S}_k\}.$$
We shall prove later that there is a subgroup $H\in {\cal S}_1$ with $\kappa _1=|H+S|-|H|.$ A maximal such a  subgroup  will be called a {\em hyper-atom}. More formal definitions will be given later.
In Section 3, we  prove the existence of hyper-atoms and obtain the following
result:

Assume that  $| S | \leq (|G|+1)/2$ and that $\kappa _2 (S)\le
|S|-1$ and let $H$ be a hyper-atom of $S$. Then $\phi _H (S)$ is either an arithmetic progression or $\kappa _2(\phi _H (S))\ge |\phi _H (S)|.$

Let  $H$ be a subgroup of an abelian group $G$. A nonempty intersection of some $H$-coset
with $A$ will be called an $H$-{\em component} of $A$.
The set of $H$-components of $A$ will be denoted by ${\cal C}_A.$
A partition of $A$ into its $H$-components
 will be called an $H$-{\em decomposition}
of $A.$ An $H$-periodic $H$-component will be called full.
A set  $A$ is said to be {\em $H$-quasi-periodic} if it has exactly one non-full $H$-component.
This component will be denoted by $A_\emptyset.$

We shall say that $A$ is
 an $H$-{\em modular-progression } if  $\phi _H(A)$ is an arithmetic progression.
 Two $H$-quasi-periodic { modular-progressions } $A$ and $B$ will be called {\em similar} if
 $\phi _H(A)$ and $\phi _H(B)$ are arithmetic progression with the same difference
 such that
 $\phi _H(A_\emptyset)$  and $\phi _H(A_\emptyset)$ are respectively  initial elements of $\phi _H(A)$ and $\phi _H(B)$.

In Section 5, we   apply the
global isoperimetric methodology introduced in \cite{hiso2007}
to prove the following Vosper type result:

Let $T$ be  a finite  subset  of
 $G$ generated by a subset $S$ such that
 $|S|\le |T|$,  $S+T$ is aperiodic,  $0\in S\cap T$ and
 $$ \frac{2|G|+2}3\ge |T+S|= |T|+|S|-1.$$
 Let $H$ be a  hyper-atom of $S$. Then $|H|\ge 2$ and moreover $T$ and $S$ are similar $H$-quasi-periodic modular progressions.


In the investigation of $T+S,$ we can assume
without loss of generality,  that $0\in T\cap S$ and $\subgp{T\cup S}=G.$

Let $T$ and $S$ be  a finite   subsets  of an abelian
 $G$ generated by $T\cup S$ such that $S$ is not an  arithmetic progression. Also, assume that $T+S$ is aperiodic, $2\le |S|\le |T|$,   $0\in T\cap S$ and
 $ |G|-2\ge |T+S|= |T|+|S|-1.$

 Kemperman's Structure Theorem states that there exists a nonzero subgroup $H$
 such that $T$ and $S$ are  $H$-quasi-periodic, $T_{\emptyset}+S_{\emptyset}$ is aperiodic and $|T_{\emptyset}+S_{\emptyset}|=|T_{\emptyset}|+|S_{\emptyset}|-1.$
 Moreover   $|\phi_H(T)+\phi_H(S)|=|\phi_H(T)|+|\phi_H(S)|-1.$

   The structure of $T$ and $S$ follow by Induction on $G/H$ and $H$, if  $\phi_H(T)+\phi_H(S)$
   is aperiodic. In order to solve the problem, Kemperman had to prove an other critical pair result
   where $T+S$ could be periodic, if there is  a unique expression  element of $T+S$.

Our $n-2$-theorem is the following:

There exists a nonzero subgroup $H$
 such that $T$ and $S$ are  $H$-quasi-periodic, $T_{\emptyset}+S_{\emptyset}$ is aperiodic and $|T_{\emptyset}+S_{\emptyset}|=|T_{\emptyset}|+|S_{\emptyset}|-1.$
 moreover one the following holds:
 \begin{itemize}
   \item[(i)]  $\phi_H(S)=\{0\}.$

    \item[(ii)]
     $\phi_H(T)+\phi _H(S)=G/H$   and
     $\phi_H(T_{\emptyset})+\phi _H(S_{\emptyset})$ is a   unique expression  element of  the factorization.
     \item[(iii)]
   $T$ and $S$ are similar $H$-quasi-periodic modular progressions.
\end{itemize}

Since each of  the conditions (i), (ii) and (ii) implies that $|T+S|= |T|+|S|-1$, our description requires no induction on $G/H.$
As one could expect, the ${n-2}$-Theorem, implies very easily Kemperman's Structure Theorem \cite{kempacta}
and its dual reconstruction given by Lev in \cite{lev}. One has just to deal with the two trivial cases $|S|=1$
and $|G\setminus (S+T)|=1.$ One needs also an easy problem that appears during the recursive procedure:
$S+T$ is periodic and contains a unique expression element.

The organization of the paper is the following:

Section 2 presents some preliminaries. In Section 3, we prove a basic property of hyper-atoms. In Section 4, we describe $T$  when
$S$ is a quasi-periodic modular progression. In Section 5, we prove  the $\frac{2n}3$-Theorem. In Section 6, we prove  the ${n-2}$-Theorem. In the last section, we investigate the strong isoperimetric property.
Since almost all the ingredients of our proofs work in if $S$ is a normal subset of a non necessarily group (for every $x,$ $xS=Sx$),
we shall investigate the strong isoperimetric property in this more general context.

\section{Terminology and preliminaries}

Recall the following result:
\begin{theirlemma}(folklore)\cite{natlivre}
Let $G$ be a finite group and let
 $A$ and $B$  be  subsets
 such that $|A|+|B|\ge |G|+t,$ where $t$ is a positive integer.
 Then  every element of $G$ has $t$ distinct representations of the form $x+y,$ where  $x\in A$ and  $x\in B$.

\label{prehistorical}
 \end{theirlemma}

The following lemma could be  known:

\begin{lemma}\label{APC}
Let $G$ be a cyclic  group generated by an element $d\in G$ and let $P\subset G$  be an   arithmetic progression  with difference $d.$
 Let $X$ be a nonempty subset of $G.$ Then $|X+P|\ge \min ( |G|,|X|+|P|-1)$.
If $|X+P|=|X|+|P|-1,$ then $X$ is an arithmetic progression with difference $d$ if one of following hold:
 \begin{itemize}
   \item[(i)] $|X+P|\le |G|-1$ or $|P|=2.$
   \item [(ii)]  For some $y\in X+P,$ $|(y-P)\cap X|=1.$
 \end{itemize}

\end{lemma}

\begin{proof}
Formulae (i)  is an easy exercise. Assume that $|X+P|=|X|+|P|-1$.

 Assume first that $|G|> |X+P|.$ Without loss of generality, we may take $P=k\{0,d\},$ where $k=|P|-1.$
 In order to have $|X+P|=|X|+k,$ we must have $|X+\{0,d\}|=|X|+1.$ Hence $X$ is an arithmetic progression with difference $d.$
 Assume now that $X+P=G.$ if $|Y|=2,$ then $|X|=|G|-1,$  and hence $X$ is an arithmetic progression with difference $d.$ Take $(y-P)\cap X=\{z\}.$ Without loss of generality, we may take $y=0.$ Put $\overline{P}=G\setminus P.$

 Clearly we have $X\subset \overline{P}.$ Since $|X|=|G|-|P|+1,$ we must have $X= \overline{P} \cup  \{z\}.$ Observe that
 $\overline{P}$ is an arithmetic progression with difference $d$. 
 The condition $(z-P)\cap \overline{P}=\emptyset$ forces that $z$ is an extremity of $P.$
It follows that $X$
 is an arithmetic progression with difference $d.$
 \end{proof}

The isoperimetric method is a global approach introduced by the author,
which
derive additive inequalities from the properties of  fragments
and atoms.
 The reader may refer to the recent paper \cite{hiso2007} for an introduction to the applications of this method.

Throughout the remaining of this section,  $G$ denotes a non-null abelian  group and  $S$ denotes a generating subset of $G$ with $0\in S.$

 For a subset $X\subset G$, we define the {\em boundary} of $X$ as $\partial _S(X)=(X+S)\setminus X.$ The boundary of $X$ with respect to
 $-S$ will be written  $\partial^- _S(X)$.
 We define the {\em co-image} of $X$ as $\nabla _S (X)=G\setminus (X+S)$.
 The co-image of $X$ with respect to
 $-S$ will be written  $\nabla^- _S(X).$ A subset $X$ with $|\nabla (X)|\ge |X|$ will be called
 {\em faithful} with respect to $S.$ The reference to $S$ could be omitted.

Notice that faithful subsets play an important role in the nonabelian case.

The next lemma is related to a notion introduced by Lee \cite{lee}:

\begin{theirlemma}\cite{balart}{Let  $X$ be a
 subset of $G$. Then $\nabla ^-(\nabla (X))+S=X+S$. \label{lee}}
\end{theirlemma}
\begin{proof}
Clearly $X\subset \nabla ^-(\nabla (X)),$ and hence $X+S\subset \nabla ^-(\nabla (X))+S$.
Put $Y= \nabla ^-(\nabla (X))\setminus X.$ One can see easily that  $(Y+S)\cap \nabla (X)=\emptyset$, and hence $Y+S\subset X+S$.\end{proof}

 We shall say that a subset
$X$ induces a {\em $k$-separation} if $ |X|\geq k$ and $|\nabla (X)|\geq
k$. We shall say that $S$ is $k$-separable if some $X$ induces a
$k$-separation.

Suppose that  $S$ is $k$-separable.
 The {\em $kth$-connectivity}
of $S$
 is defined  as
$$
\kappa _k (S )=\min  \{|\partial (X)|\   :  \ \
\infty >|X|\geq k \ {\rm and}\ |\nabla (X)|\ge k\}.
$$

Clearly
$\kappa _1(S) \le \ldots \le \kappa _k(S).$

 A finite subset $X$ of $G$ such that $|X|\ge k$,
$|\nabla (X)|\ge k$ and $|\partial (X)|=\kappa _k(S)$ is
called a {\em $k$-fragment} of $S$. A $k$-fragment with minimum
cardinality is called a {\em $k$-atom}.

It will be helpful to have in mind the  following will known lemma implicit in \cite{halgebra}:

\begin{lemma} \label{negative}{
Suppose that $S$ is $k$-separable and let $F$ be a $k$-fragment of $S.$ Then $-S$ is $k$-separable. Moreover  the following hold:

\begin{itemize}
  \item[(i)] $\kappa _k(S)=\kappa _{k}(-S).$
  \item[(ii)] If $G$ is finite, then $\nabla (F)$ is a $k$-fragment of $-S$.
  \item[(iii)] Any $k$-atom is faithful.
\end{itemize}

}\end{lemma}
\begin{proof}

Clearly, $$ \partial (X) \supset \partial ^{-} (\nabla (X)),$$ for any subset $X$ of $G.$
In particular, $-S$ is $k$-separable.  Notice that (i) follows from the definitions using the abelianity of
of the group.

We have $$\kappa _k(S)=|\partial (F)|\ge  | \partial ^{-} (\nabla (F)|\ge
\kappa _k(-S)=\kappa _k(S).$$  Thus (ii) holds.

In order to show (iii), we may assume that $G$ is finite. Let $A$ be a $k$-atom of $S$ and let $A'$ be a $k$-atom of $-S.$
If follows from the definitions that $-F$ is a $k$-fragment of $-S,$ if $F$ is a $k$-fragment of $S.$ Thus,
 $|A'|=|A|$ and $|\nabla (A)|=|G|-|A|-\kappa _k(S)=|G|-|A'|-\kappa _k(-S)= |\nabla (A')|.$
By (ii), we have $|A|=|A'|\le |\nabla ^{-} (A')|=|\nabla (A)|.$\end{proof}

Notice that (i) could not hold for infinite nonabelian groups and  that (iii) could not hold for finite nonabelian groups.

We shall say that $S$ is a {\em Vosper subset} if, for all $X\subset G$
with $|X|\ge 2$, we have $|X+S|\ge \min (|G|-1,|X|+|S|)$.





Let $S$ be a $k$-separable subset. Notice that
    $\kappa _k (S)$ is the maximal integer $j$
such that for every finite subset $X\subset G,$  with $|X|\geq k$,

\begin{equation}
|X+S|\geq \min \Big(|G|-k+1,|X|+j\Big).
\label{eqisoper0}
\end{equation}

Formulae (\ref{eqisoper0}) is an immediate consequence of the
definitions. We shall call (\ref{eqisoper0}) the {\em isoperimetric
inequality}. The reader may use the conclusion of this lemma as a
definition of $\kappa _k (S)$.

Let us point out that $S$ is $1$-separable if and only if $S\neq G.$ The following lemma, implicit in some previous papers, describes useful relations between $\kappa _1$ and $\kappa _2.$


\begin{lemma} \label{degenkappa}

Let  $ S$  be a  generating  subset of an abelian  group $G$ with
$0\in S$ and let $X$ be a subset of $G.$ The following holds.  \begin{itemize}
             \item[(i)]
If $S\ne G,$ then $|S|-1=|\partial (\{0\})|\ge \kappa _1(S).$
             \item[(ii)]
If $S$ is $2$-separable and $\kappa _2\le |S|-1,$ then $\kappa _2=\kappa_1.$
\item[(iii)]
Suppose that $S$ is $1$-separable and $\kappa _1\le |S|-2.$  Then $S$ is $2$-separable. Moreover $X$ is a $1$-fragment  (resp. $1$-atom) of $S$ if and only if $X$ is a $2$-fragment  (resp. $2$-atom) of $S$.
\end{itemize}
\end{lemma}

\begin{proof} Assume that $\kappa _2>\kappa _1$ and take a $1$-atom $A$ of $S$.
$|S|-1\ge \kappa _2>\kappa _1=|A+S|-|A|.$ It follows that $|A|\ge 2.$ Since $A$ is faithful,
we have $|\nabla (A)|\ge |A|\ge 2.$ Thus $\kappa _2\le |A+S|-|A|=\kappa _1,$ a contradiction.
The proof of (iii) is now obvious.
\end{proof}


The basic intersection theorem is the following:

\begin{theorem}\cite{halgebra,hiso2007}
Let $S$ be generating subset of an abelian group $G$ with $0\in S$.
 Let $A$ be a $k$-atom of $S$ and let
   $F$   be a   $k$-fragment of $S$ such that  $|A\cap F|\ge k$. Then
  $A\subset F.$
In particular,  distinct $k$-atoms intersect in  at most $k-1$
elements.

\label{inter2frag} \end{theorem}

The structure of $1$-atoms is the following:

\begin{proposition} \label{Cay}\cite{hejc2,hjct}

Let  $ S$  be a  generating subset of an abelian  group $G$ with
$0\in S$.  Let $H$ be a $1$-atom of $S$ with $0\in H$.
   Then
   $H$ is a subgroup.
 Moreover $\kappa _1(S)\geq \frac{|S|}{2}.$

Let $G_0$ be a group containing $G$ and let $T$ be a subset of $G_0.$ Let ${\cal V}=\{C\in {\cal C}_T=|C+S|<|H|\}$
Then
 \begin{equation}\label{olson}
|T+S|\geq (|{\cal C}_A|-|{\cal V}|)|H|+\sum _{ X\in  \cal V} |X|+|{\cal V}|\frac{|S|}{2}.
\end{equation}

\end{proposition}

\begin{proof}
Take $x\in H$. Since $x\in (H+x)\cap H$ and since $H+x$ is a
$1$-atom, we have $H+x=H$ by Theorem \ref{inter2frag}. Therefore
$H$ is a subgroup. Notice that $2|H|<|G|,$ by the definition of a $1$-atom. Since $S$ generates $G$, we have $|H+S|\ge 2|H|$,
and hence  $\kappa _1(S)=|H+S|-|H|\ge \frac{|S+H|}{2}\ge
\frac{|S|}{2}.$

Take $a_C\in C,$ for each component $C$ of $T.$ For every $C\in {\cal V},$ we have $$|C+S|=|C-a_S+S|\ge
|C-a_C|+\kappa _1\ge |C|+\frac{|S|}{2}.$$
Now (\ref{olson}) follows since $T+S=\bigcup _{C\in  {\cal C}_A } C+S$ is an $H$-decomposition.
\end{proof}

Recently, Balandraud introduced some isoperimetric objects and
proved  a strong form of Kneser's Theorem using {Proposition} \ref{Cay}.

 The next result is proved in \cite{Hejcvosp1}. The finite case is reported with almost the same proof in
 \cite{hactaa}. A short proof of this result is given in \cite{hiso2007}.

\begin{theorem} { \cite{Hejcvosp1,hactaa}\label{2atom} Let  $S$  be a finite generating
 $2$-separable subset of an abelian group $G$ with $0\in S$ and $\kappa _2 (S) \leq |S|-1$.
Let
 $H$ be a $2$-atom with $0\in H$. Then either $H$ is  a subgroup  or $|H|=2.$

\label{2atomejc} }
\end{theorem}

\begin{corollary} { [\cite{Hejcvosp1},Theorem 4.6]\label{ejcf}
Let  $S$  be a $2$-separable finite
subset of an
 abelian group $G$ such that $0\in S$, $|S|\leq (|G|+1)/2$ and $\kappa _2 (S) \leq |S|-1$.

 If  $S$ is not an arithmetic progression, then there
 is a subgroup $H$ which is a $2$-fragment of $S$.

\label{vosper}                      }
\end{corollary}

\begin{proof}

Suppose that $S$ is not an arithmetic progression and let $H$ be a  $2$-atom with $0\in H.$

 Assume first that  $\kappa
_2\leq |S|-2$ and let $K$ be a $1$-atom with $0\in A.$ By Proposition \ref{Cay}, $K$ is a subgroup.
By Lemma \ref{degenkappa}, $K$ is a $2$-fragment,    and the result holds.

 Assume now that  $$\kappa _2(S)=|S|-1.$$  In view of Theorem
\ref{2atomejc}, it is enough to consider the case $|H|=2$, say
$H=\{0,x\}$. Put $N=\subgp{x}.$

Decompose $S=S_0\cup \cdots \cup S_j$ modulo $N$, where $|S_0+H|\le
|S_1+H| \le \cdots \le |S_j+H|.$ We have $|S|+1=|H|+\kappa _2=|S+H|=\sum
\limits_{0\le i \le j}|S_i+\{0,x\}|.$

Then $|S_i|=|N|$, for all $i\ge 1$ and $S_0$ is an arithmetic progression with difference $x.$  We have $j\ge 1$, since
otherwise $S$ would be an arithmetic progression. In particular, $N$
is finite and proper.
 We have
$|N+S|<|G|$, since otherwise   $|S|\ge |G|-|N|+1\ge
\frac{|G|+2}{2},$ a contradiction.

By the definition of $\kappa _2$ and the structure of $S,$ we have  \begin{eqnarray*}
|S|-1=\kappa _2(S)&\le &|N+S|-|N|\\
&=& |S|+|N|-|S_0|-|N|\\ &\le& |S|-1, \end{eqnarray*}
and hence $N$ is a $2$-fragment.
\end{proof}

 Corollary \ref{vosper} was used to solve Lewin's
Conjecture on the  Frobenius number
 \cite{hactaa}.
Corollary \ref{vosper} coincides with [\cite{Hejcvosp1},Theorem
4.6]. A special case of this result is Theorem 6.6 of \cite{hactaa}.
As mentioned in \cite{hplagne}, there was a misprint in this last
statement. Indeed $|H| + |B| - 1$ should be replaced by $|H| + |B|$
in case (iii) of [ Theorem 6.6, \cite{hactaa}].

Alternative proofs of Corollary \ref{vosper}
 (with $|S|\leq |G|/2$
replacing $|S|\leq (|G|+1)/2$),  using Kermperman's Structure Theorem, were
obtained by Grynkiewicz in \cite{davdecomp} and Lev in
\cite{lev}. In the present paper, Corollary \ref{vosper} will be
one of the pieces leading to a generalization of Kemperman's Theorem.

Let $H$ be a subgroup of an abelian group $G$ and let $A$ and $X$ be  subsets of $G.$
An $H$-component $C$ of $X$ will be called $A$-{\em external}, if $C\cap (A+H)=\emptyset .$
Let $X\subset A$ such that $|X+H|=|H|.$ The $H$-component of $A$ {\em spanned} by $X$ is the component of $A$ containing $X.$

  We need the following consequence of
Menger's Theorem proved in \cite{hiso2007}. Notice that the condition $|\phi_H (S)|+|\phi_H (T)|\le |G/H|+1$  was omitted in \cite{hiso2007} but corrected in a another paper of the author generalizing the present work, is obviously needed. We prove in the last section a generalization of this result to the non-abelian case, valid without this restriction.
\begin{proposition} \cite{hiso2007}{
 Let $H$ be a subgroup  of an abelian group $G $.  Let $S$ and $T$ be finite subset of $G$ such that $0\in S,$
 $\kappa _1(\phi (S))= |\phi (S)|-1$ and $|\phi_H (S)|+|\phi_H (T)|\le |G/H|+1.$ Then there is a set ${\cal B}$ of
 $|\phi (S)|-1$ distinct $H$-components of $T$ and a family $\{D_C; C\in {\cal B}\}$ of $H$-components of $S$ such that the family  $\{{C+D_C}; {C\in \cal B}\}$
 span  distinct $T$-external components of $T+S.$

\label{strongip}}
\end{proposition}

We call the property given in Proposition  \ref{strongip}  the {\em
strong isoperimetric property}.

\section{Hyper-atoms }

In this section, we investigate  the new notion of a hyper-atom.
Let $S$ be a generating subset of an abelian group $G$ with $0\in S.$
Recall that $S$ is a Vosper subset if and only if  $S$ is non
$2$-separable or  $\kappa _2(S)\ge |S|,$ in view of  the { isoperimetric
inequality}, (\ref{eqisoper0}).  Assuming that $S$ is a $2$-separable Vosper subset, one may easily observe that
can be never an arithmetic progression.

\begin{lemma} { Let $S$ be a finite   generating Vosper subset of an
abelian group $G$ with $0 \in S$.   Let $X\subset G$ be a subset with $|X|\ge |S|$   and $|X+S|=|X|+|S|-1$.  Then, for every $y\in
S$, we have $|X+(S\setminus \{y\})|\ge |X|+|S|-2$. \label{vominus}}
\end{lemma}
\begin{proof}

 The result holds clearly if $S$ is an arithmetic progression (necessarily $S$ is not a $2$-separable subset in this case). So, we may assume that $|S|\ge 3.$ By  the definition of a Vosper subset, we have $|X+S|\ge |G|-1$. Assume first that $|X|=3$ and hence $|S|=3.$
The result holds unless $X+(S\setminus \{y\}=X.$ Assuming the last equality. Then  $X$ is a coset of some subgroup with order $3.$ Since $X+S$ is periodic, we must have $|X+S|\ge 6,$ a contradiction. So we may assume that $|X|\ge 4.$

Suppose that $|X+ (S\setminus \{y\})|\le |X|+|S|-3$
and take a $2$-subset   $R$ of $(X+S)\setminus ( X+(S\setminus
\{y\}))$. We have $R-y\subset X$. Also $(X\setminus (R-y))+S\subset
(X+S)\setminus R$. Thus $|(X\setminus (R-y))+S|\le |X|+|S|-3\le |G|-2$, contradicting the definition of a Vosper subset.\end{proof}

Let us prove a lemma about  {fragments in quotient groups}.

\begin{lemma} { Let $G $ be an abelian group and let $S$ be a
    finite  generating subset  $0\in S$ and $S\neq G$. Let $H$ be a
subgroup which is a $1$-fragment.  Then $H$ is faithful and
\begin{equation}\label{cosetgraph}
\kappa _1(\phi _H (S))=  |\phi _H (S)|-1.
\end{equation}
 Let $K$ be a subgroup which is a $1$-fragment of
$\phi _H (S)$ and assume that $H$ is a non-null subgroup. Then $\phi _H ^{-1}(K)$ is a $2$-fragment of $S$.
}\label{quotient}
\end{lemma}
\begin{proof}

Since $|G|>|H+S|,$ we have We have $|\nabla (H)|= |H+S|-|H|\ge |H|.$

Therefore $\phi _H (S)\ne G/H$,
 and hence $\phi _H (S)$ is
$1$-separable.
 Put $|\phi _H (S)|=u+1,$ so $\kappa _2=u|H|$.

Let $X\subset G/H$ be such that  $X+\phi _H (S)\neq G/H$. Clearly
$\phi _H^{-1} (X)+S\neq G$. Then $|\phi _H^{-1} (X)+S|\ge |\phi _H^{-1}
(X)|+\kappa _2(S)= |\phi _H^{-1} (X)|+u|H|.$

It follows that $|X+\phi _H (S)||H|\ge |X||H|+u|H|.$ Hence $\kappa
_1(\phi _H (S))\ge u=|\phi _H (S)|-1$. The reverse inequality is obvious
and follows by Lemma \ref{degenkappa}. This proves (\ref{cosetgraph}).

Let $K$ be a subgroup which is a $1$-fragment of $\phi _H (S)$. Then
$|K+\phi _H (S)|=|K|+u$. Thus $|\phi _H ^{-1} (K)+S|=|K||H|+u|H|.$ By Lemma \ref{degenkappa}, $\phi _H ^{-1} (K)$ is a $2$-fragment.\end{proof}


Let $S$ be a finite generating proper subset of an abelian group $G$ with $0 \in S.$ Proposition \ref{Cay} states that there is a $1$-atom
of $S$ which is a subgroup. A maximal subgroup which is a $1$-fragment will be called a {\em hyper-atom} of $S$. This
definition may be adapted to non-abelian groups. As we shall see, the hyper-atom is more closely related to
the critical pair theory than the $2$-atom.

\begin{theorem}\label{hyperatom}
Let $S$ be a finite $2$-separable generating subset of an abelian group $G$ such
that $0 \in S,$  $| S | \leq (|G|+1)/2$ and $\kappa _2 (S)\le
|S|-1.$ Let $H$ be a hyper-atom of $S$. Then $|H|\ge 2.$
Moreover $\phi _H (S)$ is either an arithmetic progression or a Vosper
subset.
\end{theorem}

\begin{proof}
 By Lemma \ref{degenkappa}, $\kappa _2 (S)=\kappa _1 (S).$
Let us show that \begin{equation}\label{referee}
2|\phi _H (S)|-1\le |G/H|.\end{equation}
 Clearly we may
assume that $G$ is finite.

Observe that $2|S+H|-2|H|=2\kappa _1\le 2|S|-2< |G|.$ It follows, since$|S+H|$
is a multiple of $|H|$, that $2|S+H|\le  |G|+|H|,$ and hence (\ref{referee}) holds.

 Suppose now that  $\phi _H (S)$ is not a Vosper subset. By the
 definition of a Vosper subset,  $\phi _H (S)$ is $2$-separable and $\kappa _2(\phi _H(S))\le |\phi _H(S)|-1.$

 Observe that  $\phi _H(S)$ can not
have  a $1$-fragment $M$ which is a non-zero subgroup. Otherwise  by Lemma
\ref{quotient}, $\phi _H ^{-1}(M)$ is a $2$-fragment of
$S$  strictly containing $H$, contradicting the maximality of $H$. By (\ref{referee}) and
Corollary \ref{vosper}, $\phi _H(S)$ is an arithmetic progression.\end{proof}

Theorem \ref{hyperatom} implies  a result proved by Plagne and the author \cite{hplagne}
and some extensions of it,  proved using
Kermperman's Theory,  obtained  by Grynkiewicz in
\cite{davdecomp} and Lev in \cite{lev}.

 The two main new facts
in Theorem \ref{hyperatom} are:
\begin{itemize}
\item
The subgroup $H$ in Theorem \ref{hyperatom}  is well described
as a hyper-atom.
\item
The equality $|H+S|-|H|=\kappa _1$ is   more precise
than the inequality $|H+S|\le |H|+|S|-1$ in the previous results.  This equality will be
needed later.
\end{itemize}


\section{Pairs involving a quasi-periodic modular progression}

We shall deal with sets not containing necessarily $0$.  The  important group in the isoperimetric approach is $\subgp{S-S}.$  It is easy to show that  $\subgp{S}=\subgp{S-S},$ when $S$ contains $0.$ We shall write $$T^S=(T+\subgp{S})\setminus (T+S).$$ If $0\in S$ and  $\subgp{S}=G$, then $T^S=\nabla _S (T).$

\begin{lemma}\label{nongenerating}
Let $S$ and $T$ be finite non-empty subsets of an abelian group $G$ such that $0\in T,$
$S+T$ is aperiodic and $|S+T|= |S|+|T|-1$.
  If  $T\not\subset \subgp{S-S},$
then $T$ is $\subgp{S-S}$-quasi-periodic. Moreover,   $T_\emptyset+S$ is aperiodic and $|T_\emptyset+S|=|T_\emptyset|+|S|-1$.
\end{lemma}

\begin{proof}
The case $|S|=1$ is trivial. Assume that $|S|\ge 2$ and put $M=\subgp{S-S}.$
Choose an  $a\in S$ and put $X=S-a.$ Since $S-S=X-X,$ we have $M\subset \subgp{X}.$
The other inclusion follows since $X\subset S-S.$

 Put ${\cal W}=\{C \in {\cal C}_T :  |C+X|<|M|\}.$
Since $0\in T$ and $T\not\subset M,$ we have  $|\phi _M (T)|\ge 2.$
By (\ref{olson}),

$$|T+S|= |T+X|
\ge |T|+|{\cal W}|\frac{|S|}{2}.$$
It follows that ${\cal W}=\{W\}$, form some $W\in {\cal C}_T.$
Clearly $W=T_\emptyset .$
Since $T+S$ is aperiodic, $T_\emptyset+S$ must be aperiodic. By Kneser's Theorem, $|T_\emptyset+S|\ge |T_\emptyset|+|S|-1$.
 Therefore, $$|T|+|S|-1=|T+S|\ge (\sum _{C\in {\cal C}_T \setminus \{T_\emptyset\}}|C+S|)+|T_\emptyset+S|
  \ge (\sum _{C\in {\cal C}_T \setminus \{T_\emptyset\}}|C|)+|T_\emptyset|+|S|-1\ge|T|+|S|-1.$$ The result is now obvious.
 \end{proof}

\begin{lemma}\label{transfer}
Let $S$ be an  $H$-quasi-periodic modular progression  generating an abelian group $G$ with $0\in S.$
Let  $T$ be a finite subset of  $G$ such that
$S+T$ is aperiodic and $|S|+|T|-1.$  Then $T$ is an $H$-quasi-periodic modular progression similar to $S.$
\end{lemma}

\begin{proof}

Put $|\phi _H(S)|=u$ and $|\phi _H(t)|=t.$ Take a difference $d$ of  $\phi _H(S)$ such that as $\phi _H (S_\emptyset)$ is a first element. Since $S+T$ is aperiodic, we must have $|G/H|\ge t+1+u.$

Notice that $(S\setminus S_\emptyset)+T$ is $H$-periodic and that for every component $Z$  of $S+T,$  we have $|Z|\ge \min (|T_\emptyset|,|S_\emptyset|)$.
By Lemma \ref{APC}, $|\phi_H(S+T)|\ge  t+u-1$ and $|\phi_H((S\setminus S_\emptyset)+T)|\ge  t+u-2.$ Then $T+S$ has $t+u-2$ full components,
since $(S\setminus S_\emptyset)+T$ is $H$-periodic.
We must have  $|\phi_H(S+T)|=t+u-1,$ since otherwise $T+S$ would have two more components and hence $$|T+S|\ge (t+u-2)|H|+|T_\emptyset|+|S_\emptyset|=(t-1)|H|+|T_\emptyset|+(u-1)|H|+|S_\emptyset|\ge
|T|+|S|,$$ a contradiction. Since $S+T$ is aperiodic and $(S\setminus S_\emptyset)+T$ is $H$-periodic, $\phi _H(T+S)$ must have
a unique expression  element. By Lemma \ref{APC}, $\phi _H (T)$ is a progression with  difference $d,$ having
$\phi _H (T_\emptyset)$  as first element. The possibility where $\phi _H (T_\emptyset)$  is  the last  element implies that $T+S$ is periodic.
Since $S+T$ is aperiodic, $T_\emptyset+S_\emptyset$ must be aperiodic. By Kneser's Theorem, $|T_\emptyset+S_\emptyset|\ge |T_\emptyset|+|S_\emptyset|-1.$
Now we have
  $$|S|+|T|-1=|T+S|\ge (t+u-2)|H|+|T_\emptyset+S_\emptyset|\ge (t+u-2)|H|+|S_\emptyset|+|T_\emptyset|-1\ge |S|+|T|-1.$$

In particular, $|T|= (t-1)|H|+|T_\emptyset|.$\end{proof}


\begin{lemma}\label{tpowers}
Let $S$ and $T$ be subsets of a finite  abelian group $G$, generated by $S,$ such that
 $S+T$ is aperiodic,  $0\in S\cap T$ and $|S+T|=|S|+|T|-1$. Then
 $T^S -S$ is aperiodic and  $|T^S -S|=|T^S |+|S|-1$.

 \end{lemma}

 \begin{proof} The set $T^S -S$ is aperiodic by Lemma \ref{lee}. Clearly $T^S -S\subset G\setminus T.$ Thus, $$|T^S -S|\le |G|-|T|=|G\setminus (S+T)|+|S+T|-|T|=|T^S |+|S|+|T|-1-|T|=|T^S |+|S|-1.$$ By Kneser's Theorem, we have
$|T^S -S|\ge |T^S |+|S|-1.$\end{proof}

\section{The $\frac{2n}3$-Theorem }

The following result encodes efficiently the critical pair Theory.

\begin{theorem}\label{twothird}
Let $S$ be  a finite  generating subset  of
 an
   abelian group $G$  such that $S$ is not an arithmetic progression.
Let $T$ be  a finite  subset  of
 $G$ such that
 $|S|\le |T|$,  $S+T$ is aperiodic,  $0\in S\cap T$ and
 $$ \frac{2|G|+2}3\ge |S+T|= |S|+|T|-1.$$

 Let $H$ be a hyper-atom. Then $H$ is a nonzero subgroup. Moreover $S$ and $T$ are similar $H$-quasi-periodic modular  progressions.
\end{theorem}

\begin{proof}

Set $|G|=n$, $h=|H|$,
$|\phi _H(S)|=u+1$, $|\phi _H(T)|=t+1$ and
$q=\frac{n}{h}$.

We have $|S|\le \frac{|S|+|T|}2\le \lfloor\frac{2n+5}6\rfloor< \frac{n+1}2$.
Since $S$ is not an arithmetic progression, we have $|S|\ge 3.$ Thus  $|G|>|S+T|\ge 5.$
We have $|T^S|\ge \frac{n-2}3>1,$ and hence $S$ is $2$-separable.
Thus, $\kappa _2(S)\le
|S|-1.$
 By Theorem \ref{hyperatom}, $|H|\ge 2.$

Choose   an $H$-component of $S$ with a maximal cardinality  $S_+$ and
and   an $H$-component of $S\setminus S_+$ with a minimal cardinality.
If $u\ge 2,$ we choose also an $H$-component of $S\setminus (S_+\cup S_-)$ with a minimal cardinality $S_{+-}$.
Without loss of generality, we shall assume that
 $0\in S_+$.

By the definition, we have
$u|H|=|H+S|-|H|=\kappa _1(S)\le |S|-1.$ It follows that for any subset ${\cal X}\subset {\cal C_S},$
$\sum _{C\in {\cal X}}(|H|-|C|)\le |H+S|-|H|\le |H|-1.$ Thus
\begin{equation}\label{plein}
|{\cal X}|{\max _{C\in {\cal X}} |C|} \ge\sum _{C\in {\cal X}}|C|\ge |{\cal X}||H|-(|H|-1)= (|{\cal X}|-1)|H|+1
\end{equation}

By an {\em internal} component, we shall mean an $H$-components of $T+S$ contained in $T+H$.
The set of internal components of $T$ will denoted by ${\cal I}.$
By an {\em external} component, we shall mean an $H$-component of $T+S$ disjoint from $T+H$.
Let ${\cal F}$ denotes the set of the full internal components. By ${\cal V},$ we shall denote the set of the non-full internal component. Clearly we have ${\cal I}={\cal V}\cup {\cal F}.$  By ${\cal E},$ we shall denote the set of the external components.

We shall use the following trivial observation, without any reference:

If $X\in {\cal V},$ then $|C+S_+|<|H|,$ where $C$ the component of $T$ contained in $X.$

Since $S$ generates $G$ and $H$ is proper, it follows that $u\ge 1.$
By (\ref{plein}), $|S_+|>\frac{h}2$.
Thus, $\subgp{S_+}=H$, by Lemma \ref{prehistorical}.
By (\ref{olson}), \begin{eqnarray}|T+S|
&=& \sum _{C\in {\cal F}}|C|+\sum _{C\in {\cal V}}|T+S_+|+\sum _{C\in {\cal E}}|C|
\nonumber \\&\ge&
 |{\cal F}||H|+\sum _{C\in {\cal V}}|C|+|{\cal V}|\frac{|S_+|}{2}+\sum _{C\in {\cal E}}|C|.\label{X11}
\end{eqnarray}

We have $(t+1)h\ge |T|\ge |S|> \kappa _2(S)=uh$, and hence
$$t\ge u.$$ Since  $\frac{n}3>|S|-1\ge \kappa _1(S)=|H+S|-|H|\ge h=\frac{n}q$, we must have $q\ge 4.$

{\bf Claim } 0: $ t+1+u\le q.$

Suppose the contrary. Then by Lemma \ref{prehistorical}, every element of $G/H$ has two distinct expressions. In particular,
$|E|\ge |S_{+-}|,$ for every external component $E,$ if $u\ge 2.$ Observe that any internal component $I$ contains a set of the form
$C_0+S_+,$ where $C_0$ is a component of $T$.   In particular,
$|I|\ge |S_{+}|.$

 Assume first that $u\ge 2.$
 By (\ref{plein}), $$2|S_+|\ge |S_+|+|S_{+-}|\ge \frac{2(|S_+|+|S_{+-}|+|S_-|)}{3}\ge \frac{2(2h+1)}{3}.$$
 Therefore we have
\begin{eqnarray*}
\frac{2n+2}{3}&\ge&|S+T|\\
&=& \sum _{C\in {\cal I}}|C|+\sum _{C\in {\cal E}}|C|\\
&\ge &(t+1) |S_+|+(q-t-1)|S_{+-}|\\
&=& (2t+2-q)|S_+|+(q-t-1)(|S_+|+|S_{+-}|)\\
&\ge& (2t+2-q)\frac{2h+1}{3}+2\frac{2h+1}{3}(q-t-1)
\\
&=&q\frac{2h+1}{3},\\
&\ge&\frac{2n}{3}+q/3,
\end{eqnarray*}
a contradiction, noticing that $q< t+1+u\le 2t+1.$

Assume now that  $u=1.$ We have necessarily $q= t+1$ and ${\cal E}=\emptyset$.

We must have $|{\cal V}|\le 3, $  since otherwise by (\ref{X11}), $|T+S|\ge |T|+|{\cal V}|\frac{|S_+|}2\ge |T|+|S|,$ a contradiction.
We must have $|{\cal V}|= 3, $  since otherwise by (\ref{X11}), $$|T+S|\ge (q-2)h+2|S_+|\ge(q-1)h+1=n-h+1=h(q-1)+1>
\frac{3n}4+1, $$ a contradiction.  Then $|{\cal F}|=q-|{\cal W}|\ge 4-3=1$.

Since $\subgp{S}=G$ and $u=1$, we have
$\subgp{\phi _H(S)}=\subgp{\phi _H(S_1)}=G/H$, and hence there is a component $T_0 \in {\cal F}$ such that $T_0+S_1\subset V,$  for some
 $V\in {\cal V}.$

By Lemma \ref{prehistorical}, $ |T_0| +|S_1|\le h$.  Thus by (\ref{X11}), $$|T+S|\ge
(|{\cal F}|-1)|H|+|H|+\sum _{C\in {\cal V}}|C|+\frac{3|S_+|}2\ge
|T|-|T_0|+|T_0|+|S_1|+\frac{3|S_+|}2>|T|+|S|,$$
a contradiction. The claim is proved.

{\bf  Claim} 1: $|\phi _H (S+T)|=|\phi _H (S)|+|\phi _H (T)|-1.$

By Claim 0, (\ref{cosetgraph}) and (\ref{eqisoper0}), we have
$$|\phi _H(S+T)|\ge \min(q,t+1+u)=t+1+u.$$

  By Lemma \ref{quotient},
$\kappa _1(\phi _H (S))=  |\phi _H (S)|-1$.

By Proposition \ref{strongip}, there  is a set ${\cal B}$ of
 $u$ distinct $H$-components of $T$ and a family $\{D_C; C\in {\cal B}\}$ of $H$-components of $S$ such that the family  $\{{C+D_C}; C\in {\cal B}\}$
 span  $u$ distinct $T$-external components of $T+S.$  Put  $c=|\phi _H (S+T)|-|\phi _H (S)|.$ Observe that any external component
 has a cardinality not less than $|S_-|.$

\begin{eqnarray*}
|S+T| &\ge& \sum _{C\in {\cal C}_T\setminus {\cal B}}|C+S_+|+\sum _{C\in  {\cal B}}|C+S_+|+\sum _{C\in  {\cal B}} |C+D_C|+c|S_-|\nonumber\\
&\ge& \sum _{C\in {\cal C}_T\setminus {\cal B}}|C+S_+|+\sum _{C\in  {\cal B}} |C+D_C|+\sum _{C\in  {\cal B}}|C+S_+|+c|S_-|\label{AP1}\\
&\ge& |T|+u|S_{0}|+c|S_-|. \label{AP2}
\end{eqnarray*}

We must have
 $c=0,$ since otherwise  $|T+S| \ge |T|+u|S_+|+|S_-| \ge |T|+|S|,$
a contradiction. Thus $$|\phi _H (S+T)|=|\phi _H (S)|+|\phi _H (T)|-1.$$

{\bf Claim } 2: Assume that $u\ge 2.$ Then there is at most one external component with size less than $ |S_{+-}|.$  In particular,
$\sum _{C \in {\cal E}}|C|\ge |S_{-}|+(u-1)|S_{+-}|$.

By Theorem \ref{hyperatom}, $\phi _H(S)$ is  an arithmetic progression or a Vosper subset. Let us show that
\begin{equation}
|\phi _H (T)+\phi _H (S\setminus S_-)|\ge t+u \label{omit}
\end{equation}

Observe that (\ref{omit}) is obvious if $\phi _H(S)$ is  an arithmetic progression, in view of Claim 0,
and follows  by Lemma \ref{vominus} if $\phi _H(S)$  is a Vosper subset in view of Claim 1. Claim 2 follows now.

{\bf Claim } 3: If $u\ge 2$ then $q-1\ge t+u+2$.

Assume  that $u\ge 2$ and let $T_+$ denotes an $H$-component of $T$ with a maximal cardinality.  We must have
 \begin{equation}  |{\cal F}|\ge 2.\label{EQP>1}\end{equation}
Suppose the contrary.  By (\ref{plein}), we have $|S_+|\ge \frac{2}{3}(|S_+|+|S_+|+|S_+|)\ge \frac{2h+1}{3}.$
By Lemma \ref{prehistorical},  $|C\cap T|<\frac{h}{3}<\frac{|S_+|}2$ for every $C\in {\cal V}$. By Claim 2 and (\ref{plein}), \begin{eqnarray*}
2|T|>|S+T|&\ge& \sum _{C\in {\cal V}}|(C\cap T)+S_+|+ \sum _{C\in {\cal F}}|C|+ \sum _{C\in {\cal E}}|C|\nonumber\\
&\ge& \sum _{C\in {\cal V}}2|C\cap T|+ |{\cal F}||H|+|S_{+-}|+ |S_-|\label{EQV2T0}
\\ &\ge& \sum _{C\in {\cal V}}2|C\cap T|+ (|{\cal F}|+1)h+1>2|T|,\label{EQV2T}
\end{eqnarray*}
a contradiction.
We have, using (\ref{plein}),
 \begin{equation}\label{ref2}2|S_+|\ge |S_+|+|S_{+-}|\ge \frac{2}{3}(|S_-|+|S_{+-}|+|S_{u-2}|)\ge \frac{4h+2}3.\end{equation}

Recall that  the size of an  internal component
 is not less than $|S_+|.$ The claim must hold since otherwise we have using Claim 2,
(\ref{omit}), (\ref{ref2}) and Claim 3:
  \begin{eqnarray*}
|S+T|&=&  \sum _{C \in {\cal F}} |C|+\sum _{C \in {\cal V}} |C|+\sum _{C \in {\cal E}}|C|\\
&\ge& 2|H|+(t-1)|S_+|+(u-1)|S_{+-}|+|S_-|\\
&=& 2h+ (t-2)|S_+|+(u-2)|S_{+-}|+(|S_+|+|S_{+-}|+|S_-|)\\
&=&  2h+(t-u)|S_+|+(u-2)(|S_+|+|S_{+-}|)+(|S_+|+|S_{+-}|+|S_-|)\\
&\ge& 2h+ (t-u)\frac{2h+1}3+\frac{(4h+2)(u-2)}3+2h+1\\
&\ge&(t+u+2)\frac{2h}3+1\ge q\frac{2h}3+1=\frac{2n}{3}+1,
\end{eqnarray*}
 a contradiction.


Suppose that  $\phi _H(S)$ is  an arithmetic progression, and hence $u\ge2$. By Theorem \ref{hyperatom}, $\phi _H(S)$
 a Vosper subset. By Claim 1 and Claim 3, we have
$q-2\ge |\phi _H (T+S)|=|\phi _H (S)|+|\phi _H (T)|-1,$ contradicting the definition of a Vosper subset.

Thus  $\phi _H(S)$ is  an arithmetic progression with difference $\phi _H(d),$ for some $d\in S.$
By Lemma \ref{APC} and by Claim 3,
 $\phi _H(T)$ is an arithmetic progression with  difference $\phi _H (d)$.

Now we shall order the $S_i$'s and $T_i$'s using the modular progression structure.

Take  $H$-decompositions $S=\bigcup \limits_{0\le i\le u}S_i$, $ T=\bigcup \limits_{0\le i\le t}T_i$
and an $H$-decomposition $S+T=\bigcup _{0\le i \le t+u}E_i$.
Since $\phi _H(-d)$ is also a difference of $\phi _H(S)$,  we may assume $0\in S_0$ and that
\begin{enumerate}
\item  $\phi _H(S_0), \cdots , \phi _H(S_u)$ is an arithmetic progression with difference $\phi _H(d)$ and $|S_0|\ge |S_u|$.
\item  $\phi _H(T_0), \cdots , \phi _H(T_t)$ is an arithmetic progression with difference $\phi _H(d)$.
\item $T_i+S_0\subset E_i$, for all $0\le i \le t$.
\item $T_{t}+S_{i}\subset E_{t+i}$, for all $1\le i \le u$.
\end{enumerate}
We shall put $Y=\{i \in [0,t] :  |E_i|<h\}.$
Since $|S_0|\ge |S_u|$, we have using (\ref{plein}), that $|S_0|>\frac{h}2$.
Thus $\subgp{S_0}=H$ by
Lemma \ref{prehistorical}.
By (\ref{olson}),

\begin{eqnarray}|T+S|&\ge & \sum _{0\le i \le t}|T_i+S_0|+\sum _{1\le i \le u}|T_t+S_i|\\
&\ge & |T|+|Y|\frac{|S_0|}{2}+\sum _{1\le i \le u}|T_t+S_i|.\label{Y11}
\end{eqnarray}

By (\ref{Y11}),
we have
$|T+S|\ge|T|+|Y|\frac{|S_0|}{2}+|S\setminus S_0|,$
and hence $|Y|\le 1.$

{\bf  Claim} 4: $Y=\emptyset$. Suppose the contrary. Then $Y=\{r\}$, for some $0\le r \le t$. Assume first $r<t$. By Lemma \ref{prehistorical},
$h\ge |T_r|+|S_0|$. Thus
\begin{eqnarray*}
|S+T|&\ge&|E_r|+  th+ |T_t+(S\setminus S_{0})|\\
&\ge & |S_{0}|+|T_r|+|S_{0}|+(t-1)h+|T_t|+ \sum _{1\le i\le u-1}|S_{i}|\\
&\ge& |T|+|S|-|S_u|+|S_0|\ge |S|+|T|,
\end{eqnarray*} a contradiction.
Then $r=t$. By Lemma \ref{prehistorical},
$h\ge |T_t|+|S_0|$. Also $|E_t|\ge |T_{t-1}+S_1|\ge |T_{t-1}|$. Hence
\begin{eqnarray*}
|S+T|&\ge& th+|E_t|+  |T_t+(S\setminus S_{0})|\\
&\ge & |T_t|+|S_{0}|+(t-1)h+|T_{t-1}|+ \sum _{1\le i\le u}|S_{i}|\\
&\ge& |T|+|S|,
\end{eqnarray*} a contradiction.

Let us show that  $|E_i|=h,$ for all $i\le t+u-1$.

Suppose  that  there is an  $r\le t+u-1$
with $|E_r|<h$. By Claim 4, $t+1\le r$.
 Thus since $(T_{t}+S_{r-t})\cup (T_{t-1}+S_{r-t+1})\subset E_r$, we have using (\ref{plein}), that
  $2h\ge |T_t|+|S_r|+|T_{t-1}|+|S_{r+1}|\ge |T_t|+|T_{t-1}|+h+1$, by Lemma \ref{prehistorical}.
  Thus   $|T+H|-|T|\ge 2h-(|T_t|+|T_{t-1}|)\ge h+1.$
  Now $|S+T|\ge |T+H|+|S\setminus S_0|\ge |T|+h+1+|S|-|S_0|>|T|+|S|,$ a contradiction.

Since $S+T$ is aperiodic, the set  $S_u+T_t$ is aperiodic. By Kneser's Theorem,
$|S_u+T_t|\ge |S_u|+|T_t|-1.$ Now we have
\begin{eqnarray*}
|S|+|T|-1=|S+T|&\ge& (t+u)h+|T_t+S_{u}|\\
&\ge & th+|S_{u}|+uh+|T_t|-1+ \sum _{1\le i\le u-1}|S_{i}|\\
&\ge& |T|+|S|-1
\end{eqnarray*}
Thus $|S\setminus S_u|=uh$ and $|T\setminus T_t|=th$.\end{proof}

Notice that  the subgroup in Theorem \ref{twothird} depends only one of the sets (namely $S$),
while the subgroup in Kemperman's Structure Theorem depends on $S$ and $T$.

\section{The   ${n-2}$-Theorem }

Let $G$ be an abelian group. A factorization $G=A+B$ will be called { \em singular} if there exists an $x\in G$ such that
$|A\cap (x-B)|=1.$ The element $x$ will be called a {\em unique expression} element of the factorization.

\begin{theorem}\label{n-2}

Let $T$ and $S$ be  a finite  subsets  of an abelian
 $G$ generated by $S\cup T$ such that $S+T$ is aperiodic, $2\le |S|\le |T|$,   $0\in S\cap T$ and
 $ |G|-2\ge |S+T|= |S|+|T|-1.$ If $S$ is not an arithmetic progression, then  there exists a nonzero subgroup $H$
 such that $S$ and $T$ are  $H$-quasi-periodic, $T_{\emptyset}+S_{\emptyset}$ is aperiodic and $|T_{\emptyset}+S_{\emptyset}|=|T_{\emptyset}|+|S_{\emptyset}|-1.$
 Moreover one the following holds:
 \begin{itemize}
   \item[(i)]  $\phi_H(S)=\{0\}.$

    \item[(ii)]
     $\phi_H(T)+\phi _H(S)=G/H.$ Moreover the factorization $\phi_H(T)+\phi _H(S)=G/H$ is singular and
     $\phi_H(T_{\emptyset})+\phi _H(S_{\emptyset})$ is a  unique expression  element of  the factorization.
\item[(iii)]
   $T$ and $S$ are similar $H$-quasi-periodic modular progressions.

\end{itemize}

\end{theorem}

\begin{proof}
Put $L=\subgp{S}$.
Assume that $T\not\subset  L$. Then (i) holds with $H=L$ by Lemma \ref{nongenerating}.

Form now on, we take $T\subset L,$ and hence  $S$ generates $G$.
Let us show that $G= \subgp{T}.$  Assuming the contrary. By Lemma \ref{nongenerating},    $S$ is $\subgp{T}$-quasi-periodic with necessarily two components. Thus $|S|\ge |\subgp{T}|+1\ge |S|+1,$ a contradiction.

Put $U=\subgp{T^S-T^S}.$ If  ${U}\neq G,$ we put $H=U.$ Then
 by Lemma \ref{nongenerating}, $S$ and $T$ are $H$-quasi-periodic.

 It follows that $T+S$ is $H$-quasi-periodic. Our hypothesis shows that $|H|>|T^S |\ge 2.$  We must have $T+S+H=G,$ otherwise $|T^S|>|H|.$ Thus $\phi_H(T)+\phi _H(S)=G/H.$
 Since $T+S$ is aperiodic, we must have $\phi_H(T_{\emptyset}+S_{\emptyset})$ is a unique expression  element of  the factorization.

 In this case (ii) holds.

Assume now that $U=G.$

Notice that $|S|+|T|+|T^S |=|S|+|T|+(|G|-|T+S|)=|G|+1$.
We consider the following cases:

{\bf Case} 1. $|S|\le |T^S |.$ Let $H$ denotes a hyper-atom of $S.$

Assume first that $|T|\le |T^S |.$ Then $|S|+|T|\le \frac{2(|S|+|T|+|T^S |)}3=\frac{2|G|+2}3$.   By Theorem \ref{twothird}, $T$ and $S$ are $H$-quasi-progressions with the same difference. Also $T_{\emptyset}+S$ is aperiodic and $|T_{\emptyset}+S_{\emptyset}|=|T_{\emptyset}|+|S_{\emptyset}|-1.$

Assume now  $|T|> |T^S |$,  and hence  $|S|+|T^S |\le \frac{2(|S|+|T|+|T^S |)}3=\frac{2|G|+2}3$.

By Theorem \ref{twothird}, $S$ a quasi-periodic modular $H$-progression, where $H$ is the hyper-atom of $S$. By Lemma \ref{transfer},
 $T$ a quasi-periodic modular $H$-progression similar to $S$. Also $T_{\emptyset}+S$ is aperiodic and $|T_{\emptyset}+S_{\emptyset}|=|T_{\emptyset}|+|S_{\emptyset}|-1.$
Thus (iii) holds.

 {\bf Case} 2. $U$ generates $G$ and $|S|> |T^S |.$  Let $H$ denotes a hyper-atom of $T^S-a,$
 for some $a\in T^S.$
 By Theorem \ref{twothird}, $S$ a quasi-periodic modular $H$-progression. By Lemma \ref{transfer},
 $T$ is a quasi-periodic modular $H$-progression similar to $S$. Also $T_{\emptyset}+S$ is aperiodic and $|T_{\emptyset}+S_{\emptyset}|=|T_{\emptyset}|+|S_{\emptyset}|-1.$ Thus (iii) holds.
 \end{proof}

Theorem \ref{n-2} involves some new simplifications:

\begin{itemize}
  \item  Kemperman's Structure Theorem reduces the structure of $S$ and $T$ to the of $H$-components $S_0$ and $T_0$ together with the structure of $\phi _H(S)$ and $\phi _H(T).$  Theorem \ref{n-2} does not refer to  the structure of $\phi _H(S)$ and $\phi _H(T).$
  \item Elements with a unique expression play an important role in the classical pair theory.
  Theorem \ref{n-2} avoids these elements.
  \item The quasi-period  described by Theorem \ref{n-2} i is either $\subgp{X-X}$ or a hyper-atom of some translate of
$X,$ where $X\in \{S,T^S\}.$.
\end{itemize}

\section{The strong isoperimetric property}

In this section, we shall assume some familiarity with  graphs. We shall assume also that the reader is aware of
the definition of $\kappa _1$ in the non-abelian case and its relation with the corresponding notion in Cayley
graphs. Also the notion of a component with respect to a normal subgroup may be defined as in the abelian case. These questions are explained in \cite{hiso2007}.
Possibly, the reader could restrict himself to the abelian case, where the notions are defined in the present paper.

 Let $V$ be a set and let $E \subset V\times V$.  The relation
$\Gamma = (V,E)$ will be called  a  {\em graph}.
 The elements of  $V$
will be called  {\em vertices}.
The elements of  $E$
will be called {\em arcs}.
 The graph $\Gamma$ is said to be
{\em reflexive} if $\{(x,x) : x\in V\} \subset E.$

Let $ a\in V$ and let $A\subset V$. The image of $a$ is by
definition  $$\Gamma (a)=\{x :  (a,x)\in E\}.$$ The image of $A$ is
by definition $$\Gamma (A)=\bigcup \limits_{x\in A} \Gamma (x).$$

 The {\it valency}
of $x$ is by definition $d _{\Gamma}(x)=|\Gamma (x)|$. We shall say that
$\Gamma$ is {\em locally finite} if  $d _{\Gamma} (x)$ is finite for
all $x$.

For $X\subset V$, the {\em
boundary} of $X$ is by definition
 $$\partial _{\Gamma}(X)= \Gamma (X)\setminus X.$$

 The $1$-{\em connectivity}
of $\Gamma$
 is defined  as

\begin{equation}
\kappa _1 (\Gamma )=\min  \{|\partial (X)|\   :  \ \
\infty >|X|\geq 1 \ {\rm and}\ |V\setminus  \Gamma(X)|\ge 1\}.
\label{eqcon}
\end{equation}

If $\Gamma$ is the Cayley graph defined by a generating subset $S$ of a group $G,$
then $\kappa _1(\Gamma)=\kappa _1(S).$ The reader may refer to \cite{hiso2007} for the last relation and for the definition of Cayley graphs.
 By a path from a vertex
$x$ to a vertex $y$, we shall
a finite sequence of arcs $(x,y_1), (y_1,y_2), \ldots ,(y_k,y)$. Let $\Gamma=(V,E)$ be graph.
Two paths from $x$ to $y$ are said to be openly disjoint if their intersection is $\{x,y\}.$
Recall the well known result:

\begin{theorem} ( Dirac-Menger) \label{menger}\cite{tv}

Let $\Gamma=(V,E)$ be a finite graph and let $k$ be a
nonnegative integer. Let $x,y\in V$ such that  $(x,y)\notin E$ and
$|\partial (X)|\ge k,$ for every subset $X\subset V$ with $x\in X$ and $y\notin X\cup \Gamma (X).$

Then there are
$k$ openly disjoint paths from $x$ to $y$.
\end{theorem}

\begin{proposition} ( The strong isoperimetric property)\label{SIPG}

Let $\Gamma=(V,E)$ be a locally finite graph and put $k=\kappa _1 (\Gamma ).$  Let $X\subset V$ be a finite subset such that $|X|\ge k$ and
 $|X|+k\le |V|.$ Then there are a subset $k$-subset  $C\subset X$ and an injection $f:C\rightarrow \partial (X)$ such that for every
 $c\in C,$ $(c,f(c))$ is an arc of $\Gamma.$
\end{proposition}

\begin{proof}
Take elements $\alpha$ and $\beta$ not contained in $V.$ Put $\hat{V}=X\cup \partial (X)\cup \{\alpha, \beta\}.$
We shall define a graph $\hat{\Gamma}$ on  $\hat{V}$ as follows:
\begin{itemize}
  \item $\hat{\Gamma}(\alpha)=X$ and $\hat{\Gamma}(\beta)=\emptyset.$
  \item $\hat{\Gamma}(x)=\Gamma (x),$ for every $x\in X.$
  \item $\hat{\Gamma}(x)=\beta,$ for every $x\in \partial(X).$
\end{itemize}

Take a subset $Y$ with $\alpha \in Y$ and $\beta \notin (Y\cup \hat{\Gamma}(Y)).$
It follows that $Y\cap (\partial (X)\cup \{\beta \})=\emptyset.$ Put $Y_0=Y\setminus \{\alpha\}.$
We have $$\hat{\Gamma}(Y)=\hat{\Gamma}(\alpha)\cup \hat{\Gamma}(Y_0)=X\cup \Gamma (Y_0).$$
Thus $\hat{\partial}(Y)=(X\setminus Y_0)\cup \partial (Y_0).$
It follows that $|\hat{\partial}(Y)|\ge | \partial (Y_0)|\ge \min (|V|-|X|,k)=k.$

By Menger's Theorem, there are $k$ openly disjoint paths from $\alpha$ to $\beta.$
By removing $\alpha$ and $\beta,$ we obtain $k$ disjoint paths of $\hat{\Gamma}$ from $X$ to $\partial (X).$
Take $k$ disjoint paths of $\hat{\Gamma} $ from $X$ to $\partial (X),$  with a minimal length sum. Each path consists of a single arc, since the last arc of one path is a  path of $\hat{\Gamma} $ from $X$ to $\partial (X).$
The injection $f$ is just the graph of these arcs.\end{proof}

The condition $|X|\ge k$ may be removed:
\begin{proposition} ( The strong isoperimetric property: second form)\label{SIP2}

Let $\Gamma=(V,E)$ be a locally finite graph and put $k=\kappa _1 (\Gamma ).$  Let $X\subset V$ be a finite subset such that $|X|\le k$
 and $|X|+k\le |V|.$ For every $x\in X$,  there are elements  $x_1, \ldots ,x_k\in X$ and distinct elements $y_1, \ldots ,y_k\in \partial (X)$ such that the following hold:

 \begin{itemize}
   \item  $(x_i,y_i)\in E,$ for all $1\le i \le k.$
   \item $\{x_1,x_2, \cdots ,x_{|X|}\}=X$.
   \item $x_i=x,$ for all $|X|\le i \le k$.
 \end{itemize}
\end{proposition}
This  form is not needed in the present work. So we leave the proof a exercise with a small hint:
Before applying Menger Theorem, the vertex $x$ should duplicated $k-|X|$ times.

\begin{proposition} \cite{hiso2007}{
 Let $H$ be a normal subgroup  of a multiplicative group $G $  and let $S$ and $T$ be finite subset of $G$ such  $1\in S,$
 $\kappa _1(\phi (S))= |\phi (S)|-1,$ and $|\phi_H (S)|+|\phi_H (T)|\le |G/H|+1.$ Then there is a set ${\cal B}$ of
 $|\phi (S)|-1$ distinct $H$-components of $T$ and a family $\{D_C; C\in {\cal B}\}$ of $H$-components of $S$ such that the family  $\{{CD_C}; {C\in \cal B}\}$
 span  distinct $T$-external components of $T+S.$

\label{strongip2}}
\end{proposition}
\begin{proof}
By Proposition \ref{SIPG}, there is a subset  $C$ of $\phi _H(T)$ and an injection $f:C\rightarrow \partial(\phi _H(T))$
$(c,f(c))$ is an arc. By the definition of the Cayley graph, $f(c)=c\phi(s_C),$ for some $s_c\in \phi (S).$ For each $c\in C,$
put $T_c=(\phi_H^{-1}(c))\cap T$ and $S_c=(\phi_H^{-1}(s_c))\cap S.$ The family $\{T_cS_c;c\in C\}$ satisfies the proposition.
\end{proof}

{\bf Acknowledgement}. The author is grateful to  an anonymous referee for many valuable comments on the first two drafts.

\end{document}